\documentclass[12pt,a4paper]{amsart}

\usepackage{amssymb}

\makeatletter
\renewcommand{\@begintheorem}[2]{
\rm \trivlist \item [\hskip \labelsep {\bf #2\ \ #1.}]
                                }
\makeatother

\makeatletter

\DeclareFontFamily{U}{cyr}{}
\DeclareFontShape{U}{cyr}{m}{n}{
  <5> wncyr5 <6> wncyr6 <7> wncyr7 <8> wncyr8 <9> wncyr9 <10->
wncyr10}{}
\DeclareMathAlphabet{\mathcyr}{U}{cyr}{m}{n}

\input cyracc.def

\newcommand{\ZZ}{{\bf Z}}
\newcommand{\QQ}{{\bf Q}}

\newcommand{\CC}{{\bf C}}

\newcommand{\FF}{{\bf F}}

\newcommand{\PP}{{\bf P}}

\newcommand{\et}{{\text{\'et}}}

\newcommand{\bes}{\begin{equation*}}
\newcommand{\ees}{\end{equation*}}

\newcommand{\Sym}{\mathop{\rm Sym}}

\newcommand{\tr}{\mathop{\rm tr}}
\newcommand{\Gal}{\mathop{\rm Gal}}
\newcommand{\trg}{{\mathop{\rm tr}}^{\rm Gal}}
\newcommand{\mn}{{\mathop{\rm min}}}

\title{Modularity of Maschke's octic and Calabi-Yau threefold}

\author{Matthias Sch\"utt}
\address{Institut f\"ur Algebraische Geometrie, Leibniz Universit\"at
  Hannover, Welfengarten 1, 30167 Hannover, Germany}
\email{schuett@math.uni-hannover.de}

\date{December 30, 2011}

\begin{document}

\begin{abstract}
We prove the modularity of Maschke's octic and two Calabi-Yau threefolds
derived from it as double octic and quotient thereof by a suitable Heisenberg group,
as conjectured by Bini and van Geemen.
The proofs rely on automorphisms of the varieties
and isogenies of K3 surfaces.
\end{abstract}

\maketitle

\section{Introduction}

In their paper \cite{BvG}, Bini and van Geemen study Maschke's octic surface $S$ and two Calabi-Yau threefolds $X$ and $Y$ which arise as double covering of $\PP^3$ branched along $S$
and quotient thereof by a suitable Heisenberg group.
Thanks to the huge automorphism group of $S$,
Bini and van Geemen are able to decompose the cohomology of either variety into small pieces, each of dimension 1, 2 or 3.
The decomposition works over $\QQ$ (for $S$ and $Y$) or $\QQ(i)$ where $i^2=-1$ (for $X$).
This leads to the problem of the corresponding Galois representations. 
The aim of this note is to prove their modularity, as conjectured by 
Bini and van Geemen.

Throughout the paper we employ the notation from \cite{BvG},
see \ref{s:not} for a brief account.
Particularly we are concerned with the following modular forms 
with notation from \cite{MFIV} indicating the level:
\[
f{120} \;\;\; \text{weight } 4;\;\;\; f{15C}, f{24b}, f{120E} \;\;\; \text{weight } 2
\]
(see \cite{MFIV} or \cite[7.2]{BvG} for Fourier coefficients).
The associated compatible systems of $\ell$-adic Galois representations 
are denoted by $V_{f,\ell}$.
 In detail we prove the following:

\subsection{Theorem}
\label{thmX}
{\sl Maschke's double octic $X$ is modular over $\QQ$. 
 Its third cohomological Galois representation decomposes as
\begin{eqnarray*}
H^3_\et(X,\QQ_\ell)  & \cong & V_{f{120},\ell} \oplus
V_{f{120E,\ell}}(-1)^{32}
\oplus
(V_{f{120E,\ell}}(-1)\otimes\chi_{-1})^{18}
\oplus V_{f{24B,\ell}}(-1)^{36}\\
&& \;\;
\oplus (V_{f{24B,\ell}}(-1)\otimes\chi_{-1})^{18}
\oplus V_{f{15C,\ell}}(-1)^{27}
\oplus (V_{f{15C,\ell}}(-1)\otimes\chi_{-1})^{18}
\end{eqnarray*}
where $\chi_{-1}$ denotes the quadratic Dirichlet character of conductor $4$.}

\subsection{Theorem}
\label{thmY}
{\sl The Calabi-Yau threefold $Y$ is modular over $\QQ$. 
Its third cohomological Galois representation decomposes as
\[
H^3_\et(Y,\QQ_\ell) \cong V_{f{120},\ell} \oplus
V_{f{120E,\ell}}(-1)^5
\oplus V_{f{24B,\ell}}(-1)^9.
\]}

\subsection{Theorem}
\label{thmS}
{\sl Maschke's octic $S$ is modular over $\QQ$. 
Its transcendental Galois representation decomposes as
\[
T_{S,\ell} \cong 
(V_{f{15},\ell})^5 \oplus
({\Sym}^2 V_{f{1200},\ell}\otimes\chi_{-15})^{18}
 \oplus
({\Sym}^2 V_{f{1200},\ell}\otimes\chi_{15})^{12}
\]
for the weight 3 form $f{15}$ of level $15$ from \cite[Table 2]{S-CM} 
and the newform $f{1200} \in S_2(\Gamma_0(1200),\chi_3)$ specified in \ref{ss:1200}.
The quadratic characters $\chi_3, \chi_{15}, \chi_{-15}$ are Legendre symbols
as indicated in the subscript.
}

\medskip

The proof of Theorem \ref{thmY} is given in section \ref{s:Y}.
It builds on effective (and feasible) explicit methods for two-dimensional Galois representations,
combined with the use of involutions on $Y$.
Section \ref{s:X} proceeds with the proof of Theorem \ref{thmX},
using similar techniques.
A notable difference is that
we first have to prove that the cohomological Galois representation decomposes
into 2-dimensional pieces over $\QQ$.
Here again automorphisms of $X$ play a crucial role.
The note concludes with the proof of Theorem \ref{thmS} in section \ref{s:S}.
It proceeds by relating $S$ to a product of an elliptic $\QQ$-curve $E$ and one of its conjugates 
through various quotient maps of K3 surfaces and eventually a Shioda-Inose structure
as elaborated in \cite{GS}.
Then this $\QQ$-curve is modular with associated newform $f{1200}$, 
and we verify the claimed equality.
As a by-product we also deduce the modularity of the K3 surfaces involved in this construction,
and of another K3 surface related to $E$ by \cite{GS}
(see Remark \ref{ss:rem}).

\subsection{Notation and Set-up}
\label{s:not}

Maschke's octic is defined by the homogeneous equation
$$
S=\left\{\sum_{i=0}^3 x_i^8\,+\,14\sum_{i<j}x_i^4x_j^4\,+\,168x_0^2x_1^2x_2^2x_3^2=0\right\}\subset\PP^3.
$$
In \cite{BvG} it is checked using Groebner basis that $S$ is smooth outside characteristics $2,3,5$, in particular over $\CC$.
$S$ comes with a large group of automorphisms $\bar G$ of order $11520$ 
which can be used to analyze its (middle) cohomology.
Bini and van Geemen compute the Picard number $\rho(S)=202$
and decompose the 100-dimensional transcendental part of $H^2(S)$ into 2- and 3-dimensional pieces over $\QQ$ (as Galois representations or Hodge structures).
Their modularity is the subject of Theorem \ref{thmS}, see Section \ref{s:S} for the proof.

The automorphism group lifts to the Calabi-Yau threefold $X$
which arises as double cover of $\PP^3$ branched along $S$:
\[
X \stackrel{2:1}{\longrightarrow} \PP^3.
\]
This double octic has quite large middle cohomology, $b_3(X)=300$,
but the automorphisms can be used to split it up completely into 2-dimensional pieces
over $\QQ(i)$.
In Section \ref{s:X} we show that 
these 2-dimensional Galois representations indeed descend to $\QQ$,
and prove their modularity (Theorem \ref{thmX}).

A non-trivial portion of $H^3(X)$ can be handled by considering another Calabi-Yau threefold $Y$ which is the desingularisation of the 
quotient of $X$ by a suitable Heisenberg group $H$ (see \cite{BvG}).
Here one has $b_3(Y)=30$, 
and $H^3(Y)$ is decomposed into 2-dimensional Galois representations
properly over $\QQ$ in \cite{BvG},
so we can attack their modularity right away. 
This is achieved in the next section, proving Theorem \ref{thmY}.

\section{Proof of Theorem \ref{thmY}}
\label{s:Y}

In this section we prove the modularity of the Calabi-Yau threefold $Y$ over $\QQ$.
Recall that $Y$ is the quotient of Maschke's double octic $X$ by a suitable Heisenberg group $H$.
In \cite{BvG} it is proved that the  Galois representation of the middle cohomology of $Y$ decomposes as
\begin{eqnarray}
\label{eq:Y}
H^3_\et(Y,\QQ_\ell) =
W_1 \oplus W_5^5 \oplus W_9^9
\end{eqnarray}
for two-dimensional Galois representations $W_1, W_5, W_9$ over $\QQ$.
These Galois representations can be determined by computing sufficiently many traces of Frobenius elements (up to $p=73$, see \ref{ss:Gal}).
In general these traces can be obtained from $\#Y(\FF_q)$ for $q=p, p^2, p^3$ by 
Lefschetz' fixed point formula.
However, since $Y$ is not a hypersurface, already counting points over $\FF_{p^2}$
soon becomes infeasible.
To circumvent this, we use involutions of $Y$ to derive all required traces from point counts over $\FF_p$ for $Y$ and two quadratic twists. 
Before going into the details we prove an auxiliary property concerning Tate twists.

\subsection{Tate twist}
\label{ss:Tate}

We start by proving that $W_5, W_9$
are Tate twists of integral Galois representations of weight 1
(i.e.~the eigenvalues have absolute value the square root of the norm).
For this purpose we invoke the fact
that the Newton polygon lies above the Hodge polygon (cf.~\cite{Mazur}).
Since we work with $H^3(Y)$,
the Hodge polygon has slopes $0, 14 \times 1, 14 \times 2, 3$.
Presently the decisive property are the multiplicities (being greater than 1)
combined with the following easy fact about 2-dimensional Galois representations:

\subsection{Fact}
\label{fact}
A two-dimensional Galois representation $Gal(\bar\QQ/\QQ)\to$ GL$_2(\QQ_\ell)$ of weight $k$ can only have Newton slopes
$0,1,\hdots,k$ and $k/2$ (if this is not already an integer).

\begin{proof}
Since the characteristic polynomial at $p$ splits over some quadratic extension $K$ of $\QQ$,
the slopes are half-integers.
In fact they are integers unless $p$ ramifies in $K$
which immediately gives slope $k/2$.
\end{proof}

\subsection{Lemma}
\label{lem:Tate}
{\sl 
$W_5, W_9$ are Tate twists of integral weight one Galois representations $U_5, U_9$:
\[
W_5 = U_5(-1),\;\;\; W_9 = U_9(-1).
\]}

\begin{proof}
Fact \ref{fact} a priori leaves
 the possibilities $0,1,3/2,2,3$ for the slopes of $W_5, W_9$ at a fixed prime $p$.
Assuming that one Newton slope were $0$, this would count with multiplicity 5 resp.~9 
in $H^3_\et(Y,\QQ_\ell)$.
In either case, this would cause the Newton polygon of $H^3(Y)$ to sit underneath the Hodge polygon.
This gives the required contradiction.
By symmetry the Newton slopes thus can only be $1,3/2,2$ at any prime $p$.
This shows that $W_5$ and $W_9$ are Tate twists.
\end{proof}

\subsection{Traces}

Since Lemma \ref{lem:Tate} limits the possible traces for $W_5$ and $W_9$ 
(using the bounds imposed by the Weil conjectures),
this in fact often enables us to compute the traces from substantially less information than
point counts up to $\FF_{p^3}$. 
Namely counting points over $\FF_p$ and $\FF_{p^2}$ we obtain a linear and a quadratic relation for the traces:
\begin{eqnarray*}
1+p+p^2+p^3 - \# Y(\FF_p)  & = & tr_p(W_1) + 5p\, tr_p(U_5) + 9p\, tr_p(U_9),\\
1+p^2+p^4+p^6 - \# Y(\FF_{p^2})  & = & tr_p(W_1)^2 + 5p^2\, tr_p(U_5)^2 + 9p^2\, tr_p(U_9)^2 - 30 p^3.
\end{eqnarray*}
For $p=7,11,13$ we verified that this gives a unique integer solution respecting the Weil bounds.
This agrees exactly with the values from \cite[7.2]{BvG}:

{\renewcommand{\arraystretch}{1.2}
$$
\begin{array}{|r||r|r|r|}\hline
p&7&11&13\\
\hline
\hline
tr_p(W_1) &0&4&54\\
\hline
tr_p(U_5)&0&-4&6\\
\hline
tr_p(U_9) &0 &4&-2\\
\hline
\end{array}
$$
}

\subsection{Involutions}

In order to compute these traces for larger $p$ without counting points over $\FF_{p^2}$
we consider the involutions 
\[
\imath_1: 
\begin{cases}
x_0\mapsto -x_0\\
x_1\mapsto -x_1
\end{cases}\;\;\; \imath_2: x_4 \mapsto -x_4.
\]
We want to compute the traces of $\imath_1^*, \imath_2^*$ on our representations $W_1, W_5, W_9$ inside $H^3_\et(Y,\QQ_\ell)$.
This could certainly be done abstractly, but we will use the composition with Frobenius at the above primes and Lefschetz' fixed point formula again.
Since the involutions and Frobenius commute, these traces are directly related to the number of fixed points  of their composition.
Presently the fixed points of $F_p\circ\imath_j$ are defined over $\FF_{p^2}$,
but they are in 1:1-correspondence with $F_p$-fixed points of a quadratic twist of $Y$ (at $(x_0, x_1)$ resp.~$x_4$) and thus easy to compute.
At the above primes we find:

{\renewcommand{\arraystretch}{1.2}
$$
\begin{array}{|r||r|r|r|}\hline
p&7&11&13\\
\hline
\hline
tr (F_p^*\circ\imath_1^*) &0&180&-102\\
\hline
tr (F_p^*\circ\imath_2^*) &0&-4&-210\\
\hline
\end{array}
$$
}

By Lemma \ref{lem:Tate} we have
\[
tr (F_p^*\circ\imath_j^*) \equiv \pm\, tr_p(W_1) \mod p \;\;\; (j=1,2).
\]
The trace at $p=11$ or $13$ thus gives the sign of $\imath_j^*$ on $W_1$.
Then the traces of $F_p^*\circ\imath_j^*$ for $p=11,13$ give a linear system
in the two unknown traces of $\imath_j^*$ on $W_5$ and $W_9$.
The unique solution is (in the order $(W_1, W_5, W_9)$)
\begin{eqnarray}
\label{eq:inv}
tr\, \imath_1^* = (1,-1,3),\;\;\; tr\, \imath_2^* = (-1,-3,-3).
\end{eqnarray}

\subsection{Traces and Fourier coefficients}
\label{ss:tr}

Thanks to the involutions $\imath_1, \imath_2$ and their traces \eqref{eq:inv},
we can compute the traces of Frobenius on $W_1, U_5, U_9$ from the number of points of $Y$ and its two quadratic twists over $\FF_p$.
For $p$ from $7$ up to $73$ we find exactly the Fourier coefficients of the newforms from Theorem \ref{thmY}.

\subsection{Galois representations: even traces}
\label{ss:even}

We are now in the position to prove that the Galois representations are isomorphic up to semi-simplification. 
There is an effective method due to Faltings-Serre-Livn\'e
that applies to $2$-adic Galois representations with even trace \cite{Livne}.
In order to see that this is applicable, we use a standard argument
that we only sketch below.

Note that any Frobenius element $F_p$ with odd trace
would imply that the kernel of the mod $2$-reduction of the Galois representation
would correspond to a Galois extension of $\QQ$ with Galois group $C_3$ or $S_3$.
Here $F_p$ represents the conjugacy class of an element of order $3$. 
Going through the list of all such extension of $\QQ$ which are unramified outside the set
of bad primes $\{2,3,5\}$ \cite{Jones},
we find that each has an order 3 Galois element represented by $F_p$ for $p\in\{7,11,13,19\}$.
Thus the fact that the traces are even at these four primes implies that the Galois representations have already even trace.

\subsection{Galois representations: isomorphic semi-simplifications}
\label{ss:Gal}

In case of even traces (and equal determinant as assured presently) the theorem of Faltings-Serre-Livn\'e states
that it suffices to compare traces at a non-cubic set of primes
for the Galois group of the compositum of all quadratic extensions of $\QQ$ which are unramified outside the bad primes \cite[Thm.~4.3]{Livne}.
Presently this extension is $K=\QQ(i,\sqrt{2},\sqrt{-3},\sqrt{5})$. 
The Galois group $Gal(K/\QQ)$ is represented by the Frobenius elements at the primes
\begin{eqnarray}
\label{eq:nc}
\{7,11,13,17,19,23,29,31,41,43,53,61,71,73,83,241\}.
\end{eqnarray}
In the non-cubic set we can omit the primes $83, 241$ from the full set of representatives
(see \cite[Prop.~4.11]{Livne}).
It follows that equality of the traces at primes up to $p=73$ suffices to deduce 
that the Galois representations in question are isomorphic (cf.~\ref{ss:tr}).
Here we use that 2-dimensional integral Galois representations of odd weight are 
automatically semi-simple.
This completes the proof of Theorem \ref{thmY}.
\qed

\section{Proof of Theorem \ref{thmX}}
\label{s:X}

With the modularity of $Y$ at our disposal,
we continue to prove the modularity of the double octic $X$, branched along Maschke's octic $S$.
Recall  from \cite{BvG} that the automorphism group of $S$, interpreted on $X$,
splits the third cohomological Galois representation of $X$ over $\QQ(i)$ as follows:
\begin{eqnarray}
\label{eq:H^3(X)}
H^3_\et(X,\QQ_\ell) &\cong&
H^3_\et(Y,\QQ_\ell) \oplus
W_{15,\ell}^{15}\,\oplus\,W_{30,\ell}^{30}\,
\oplus\;
W_{45,\ell}^{45}\;\oplus\;
(W_{45,\ell}')^{45}.
\end{eqnarray}

\subsection{Tate twist}

As in \ref{ss:Tate} we find that each summand $W$ in \eqref{eq:H^3(X)} is a Tate twist of an integral Galois representation of weight 1.
The argument a priori applies to the 2-dimensional Galois representations over $\QQ(i)$, but naturally descends to the possibly 4-dimensional Galois representations over $\QQ$.

\subsection{Traces at primes $p\equiv 1 \mod 4$}
\label{ss:aut}

Using automorphisms of $X$,
we can compute the traces on the Galois representations from \eqref{eq:H^3(X)}
at any F$_q$ for $q\equiv 1\mod 4$ from point counts over $\FF_q$
on $X$ and different twists.
In detail one can use the following automorphisms and their traces on the single representations $V_i$:

{\renewcommand{\arraystretch}{1.2}
$$
\begin{array}{|c|c||r|r|r|r|r|r|r|}
\hline
\multicolumn{2}{|c||}{} & \multicolumn{7}{|c|}{\tr \imath(V_i)}\\ 
\hline
\multicolumn{2}{|c||}{\imath} & V_1 & V_5 & V_9 & V_{15} & V_{30} & V_{45} & V_{45}'\\
\hline
\hline
\imath_1 & (x_2, x_3) \mapsto (-x_2,-x_3) & 1 & 5 & 9 & -1 & -2 & -3 & -3\\
\hline
\imath_2 & (x_1, x_2, x_3) \mapsto (-x_1, x_3, x_2) & 1 & 1 & 1 & 3 & 2 & -3 & 1\\
\hline
\imath_3 & (x_1, x_2, x_3) \mapsto (ix_1, -ix_2, -x_3) & -1 & 1 & -3 & -3 & -2 & 1 & 5\\
\hline
\end{array}
$$
}

Then compose these automorphisms with F$_q$.
The fixed point sets can be derived from $X(\FF_q)$ and the $\FF_q$-rational points
of specific quadratic or quartic twists of $X$.
Through Lefschetz' fixed point formula,
these quantities relate to the trace on $H^3_\et(X,\QQ_\ell)$.
Since F$_q$ and the automorphisms commute by choice of $q$,
this trace can be computed as
\[
\mbox{tr}\, (F_q\circ\imath)^* (H^3_\et(X,\QQ_\ell)) = \sum_i \mbox{tr}\, F_q^*(W_{i,\ell}) \mbox{tr}\, \imath^*(V_i).
\]
Note that the $4\times4$-matrix comprising the traces of the automorphisms (including the identity)
on the $V_i (i=15,30,45,45')$ is invertible of determinant $960$.
Hence we obtain the traces of F$_q^*$ on the $W_i$ from
the known Galois representation of $Y$
and  point counts over $\FF_q$
on $X$ and the three twists of $X$ as sketched above.

In particular, for all $p\equiv 1 \mod 4$ up to $p=97$ we can easily verify 
that 
\begin{itemize}
\item
$W_{15,\ell}$ and $W_{30,\ell}$ have the same trace (matching the modular form $f15C$),
\item
$W_{9,\ell}$ and $W_{45,\ell}$ have the same trace (matching the modular form $f24B$),
\item
$W_{5,\ell}$ and $W_{45,\ell}'$ have the same trace (matching the modular form $f120E$).
\end{itemize}
Note that this agrees perfectly with Theorem \ref{thmX}.
It remains to deal with (a few) primes $p\equiv 3\mod 4$.

\subsection{2-dimensional Galois representations}

We shall now prove that all the $W_i$ actually descend to 2-dimensional Galois representations over $\QQ$.
Thus we have to exclude that some form an irreducible Galois representation of dimension 4 over $\QQ$.
Note that dimension 4 would imply that Frobenius eigenvalues are duplicated at $q\equiv 1\mod 4$
and appear with opposite signs at $q\equiv 3 \mod 4$
(so that the 4-dimensional representation would have zero trace
as it equals its twist by $\chi_{-1}$).

First we use the duplication observation.
Note that  none of the above three cases have the same traces at both $p=13$ and $29$ (cf.~\cite[7.2]{BvG}).
Hence there cannot be two representations of multiplicity 45 and dimension two over $\QQ(i)$,
forming the restriction of a single 4-dimensional Galois representation over $\QQ$.
Hence, for multiplicities 15 and 45, we deduce from the odd parity of the multiplicity over $\QQ(i)$ that the representations do in fact come from 2-dimensional Galois representations over $\QQ$.
Thus only $W_{30,\ell}$ may a priori stem from a 4-dimensional Galois representation $W_\ell$ over $\QQ$.
In order to establish a contradiction, we assume that it actually does:

\subsection{Assumption}
\label{ss:ass}

{\sl
$W_{30,\ell}$ comes from an irreducible  4-dimensional Galois representation $W_\ell$ over $\QQ$.
}

\smallskip

A trivial, but crucial fact now is that for $p\equiv 3\mod 4$ not only F$_p^*$ has zero trace on $W_\ell$,
but also its composition with any automorphism of $X$.
We will use this to derive a contradiction from a detailed trace analysis at a few primes.


\subsection{Galois-trace}

The 2-dimensional Galois representations $W_i$ ($=W_{15}, W_{45}, W_{45}'$) over $\QQ(i)$ can descend
in two ways to $\QQ$, distinguished by the quadratic twist by $\chi_{-1}$.
In accordance, we split the representations $V_i$ as
\[
V_i = V_i' \oplus V_i'',\;\;\; \dim(V_i')+\dim(V_i'')=15 \;\text{~resp.~}\; 45
\]
such that the corresponding Galois representations over $\QQ$ fulfill
\[
(W_i')_\ell =W_{i,\ell,\QQ}, \;\;\; (W_i'')_\ell =  W_{i,\ell,\QQ}\otimes\chi_{-1}
\]
for a fixed 2-dimensional Galois representation $W_{i,\ell,\QQ}$ (unique up to twist by $\chi_{-1}$).
These quadratic twists interfere with the action of the automorphisms from \ref{ss:aut}.
Our first goal is to compute the trace of the automorphisms weighted by the twists;
for simplicity we shall call this the Galois-trace of the automorphism $\imath$:
\[
 \trg\, \imath(V_i) = \tr \, \imath^* (V_i') - \tr \, \imath^* (V_i'').
\]
Here is the motivation of the definition:
for any prime $p\equiv 3\mod 4$ we have
\begin{eqnarray}
\label{eq:trg}
\tr (F_p\circ\imath)^* (V_i\otimes W_i)_\ell = (\trg\, \imath(V_i)) (\tr F_p^*(W_{i,\ell,\QQ})).
\end{eqnarray}
Recall that the sum over expressions on the left hand-side can be computed from the Lefschetz' fixed point formula 
(thus point counting on a twist over $\FF_p$)
and the knowledge of the action on $H^3_\et(Y,\QQ_\ell)$.
We shall now compute the Galois traces, still assuming \ref{ss:ass},
and derive a contradiction.

\subsection{Trace analysis for small primes}
\label{ss:small}

Consider the primes $p=11,19,23,31$ which are all congruent to $3 \mod 4$.
Let $p$ be one of them and $q=p^2$.
Using the automorphisms from \ref{ss:aut} composed with F$_q$,
we can compute the traces of F$_q^*$ on the $W_{i,\ell}$ from 4 point counts over $\FF_q$.
We verify that the traces satisfy the  natural relation with the coefficients $A_p$ of the modular forms as specified above:
\[
\mbox{tr}\, F_q^*(W_{i,\ell}) =  p^2(A_p^2-2p).
\]
Conversely this gives us the possible traces of F$_p^*$ on the $W_{i,\ell}$ up to sign.
Once and for all and without loss of generality, we choose the signs at $p=11$ matching the modular forms as in \ref{ss:aut}.

\subsection{Involution $\imath_2$}
\label{ss:inv2}

For any sign choice of $\tr F_p^*(W_{i,\ell,\QQ})$ for $p=19,23,31$,
summation over \eqref{eq:trg} yields a linear system of equations with Galois traces as unknowns.
With signs chosen to match the modular forms,
one can easily check that the only integral solution respecting the trivial bounds
\[
|\trg\,\imath_2(V_i)| \leq \dim(V_i)
\]
occurs for 
\[
\trg \imath_2(V_{15}) = \trg \imath_2(V_{45}) = 1,\;\;\; \trg \imath_2(V_{45}') = -3.
\]
In fact, we also find unique traces of F$_p$ 
on the 2-dimensional Galois representations $W_{i,\ell,\QQ}$.
As desired these match the modular forms as in \ref{ss:aut}.
With the traces of F$_p$ at hand, we can compute all other Galois traces
from the corresponding system of linear equations (after \eqref{eq:trg}):

{\renewcommand{\arraystretch}{1.2}
$$
\begin{array}{|c||r|r|r|}
\hline
 & \multicolumn{3}{|c|}{\trg \imath(V_i)}\\ 
\hline
{\imath} & V_{15}  & V_{45} & V_{45}'\\
\hline
\hline
id & 9 & 9 & 9\\
\hline
\imath_1   & 9 & 9 & 9\\
\hline
\imath_2 
& 1 & 1 & -3\\
\hline
\imath_3 
& -1  & 5 & 1\\
\hline
\end{array}
$$
}

\subsection{Contradiction to Assumption \ref{ss:ass}}
\label{ss:sign}

We shall use the fact that $id$ and $\imath_1$ have the same Galois traces on $V_{15}$ to establish a contradiction.
From $\trg\, id(V_{15})=9$,
we directly infer that
\begin{eqnarray}
\label{eq:dims}
\dim V_{15}'=12, \;\;\; \dim V_{15}''=3.
\end{eqnarray}
By \ref{ss:aut} we have 
\[
\tr \,\imath_1(V_{15}) = -1,
\]
i.e.~$\imath_1^*$ has 7 eigenvalues $+1$ and 8 eigenvalues $-1$ on the total representation $V_{15}$.
Combining these eigenvalues with those of the character $\chi_{-1}$ reflected in \eqref{eq:dims},
we find at least 5 products evaluating as $-1$.
That is, $\trg\imath_1(V_{15})\leq 10\cdot(+1) + 5\cdot(-1)=5<9$, contradiction. 
Thus our original assumption \ref{ss:ass} is seen to be false, and we have proved:

\subsection{Proposition}

{\sl $H^3_\et(X,\QQ_\ell)$ decomposes into 2-dimensional Galois representations over $\QQ$.}

\medskip

The alert reader might wonder why we nonetheless have found the right traces in the above argument.
The reason is that the representations of multiplicities 15 and 30 do in fact contribute
many zero traces to the cohomology (at primes $p\equiv 3\mod 4$),
but this is not fully supported on the representation of multiplicity 30
when we take automorphisms into account.

\subsection{Galois-traces revisited}
\label{ss:rev}

We reconsider the Galois traces of the involutions,
now also taking the representation $V_{30}$ into account.
At the primes $p=11,19,23,31$, we find the same traces (up to sign) 
for $W_{15,\ell}$ and $W_{30,\ell}$.
Thus we can compute the combined Galois trace, assuming without loss of generality that 
$W_{15,\ell,\QQ}$ and $W_{30,\ell,\QQ}$ take the same sign.
In fact, the argument using $\imath_2$ in \ref{ss:inv2}
goes through just the same with the dimension bound for $V_{15}$ replaced by
\[
|\trg\,\imath_2(V_{15} \oplus V_{30})| \leq 45.
\]
In particular, we find the same Galois trace as before, but now this is supported on the 
larger representation $V_{15} \oplus V_{30}$ (so that no elementary sign argument would cause a contradiction).

\subsection{Continuation of proof}

Eventually we want to apply the Faltings-Serre-Livn\'e method from \ref{ss:Gal} again.
In view of \ref{ss:aut} and \ref{ss:rev},
only the traces at $p=43,71$ are missing from the non-cubic set  \eqref{eq:nc}.
One way to proceed would be to first verify 
that $W_{15,\ell}$ and $W_{30,\ell}$ have the same trace (up to sign)
at these two primes and then proceed as above.
This could be achieved by computing the full set of traces of F$_{p^2}$ as in \ref{ss:aut}, \ref{ss:small}.
However, since these computations tend to become fairly involved, we decided to pursue another 
line of argument that we hope to be of independent interest.
Our ideas rely on elementary linear algebra, the Weil bounds for the traces of Frobenius, and the following divisibility property at the primes congruent to $3$ modulo $4$:

\subsection{Lemma}
\label{lem:4}

{\sl The Galois representations $W_{i,\ell}$ have 4-divisible traces at primes $p\equiv 3\mod 4, p>3$.}

In view of the traces computed thus far,
the lemma will follow immediately from the following general criterion:

\subsection{Theorem}
\label{prop:4}

{\sl
Let $\rho:  \Gal(\bar\QQ/\QQ) \to \mbox{GL}_2(\QQ_\ell)$
denote an integral Galois representation of odd weight 
which is unramified outside $\{2,3,5\}$.
Assume that
\[
\tr \rho(\mbox{Fr}_{13})\equiv 0\mod 2, \;\;\; \tr \rho(\mbox{Fr}_p)\equiv 0\mod 4 \;\;\; (p=7,11,19).
\]
Then $\rho$ has even trace and moreover
\[
\tr \rho(\mbox{Fr}_p)\equiv 0\mod 4 \;\;\;\;\; \forall\, p\equiv 3\mod 4, p>3.
\] 
}

\begin{proof}
By \ref{ss:even} $\rho$ has even trace.
In particular the mod $2$-reduction $\rho_2$ has image consisting of 1 or 2 elements.
Now consider the lifts to GL$(\ZZ/4\ZZ)$, i.e.~the possible images $G$ of the mod $4$-reduction $\rho_4$.
One easily checks that $\#G\leq 32$ with every element of order $1, 2$ or $4$.
Integrality and odd weight $k$ imply that
\[
\det \rho = \chi_\ell^k
\]
for the $\ell$-adic cyclotomic character $\chi_\ell$.
In particular, we find that
\[
\det \rho_4(\mbox{Fr}_p) = p
\]
specialising to $-1$ in $\ZZ/4\ZZ$ for all $p\equiv 3\mod 4$.
A case by case-analysis teaches us that the trace of all determinant $-1$ elements in $G$ is
\begin{itemize}
\item
$0$ for lifts of the identity;
\item
$0$  for lifts of other elements in GL$(\ZZ/2\ZZ)$ iff the order is $2$;
\item
$2$ for lifts of other elements in GL$(\ZZ/2\ZZ)$ iff the order is $4$.
\end{itemize}
Hence if we assume that there is a prime $p\equiv 3\mod 4$ 
with 
$\tr \rho(\mbox{Fr}_p)\not\equiv 0\mod 4$,
then this corresponds to a Galois extension $K/\QQ$ 
unramified outside $\{2,3,5\}$ 
with Galois group $C_4$ or $D_4$
where F$_p$ induces an element of order $4$.
Going through the list of all such extensions from \cite{Jones}
we find that
\begin{itemize}
\item
either $p=7,11$ or $19$ gives a Galois element of order $4$
\item
or the extension $K/\QQ$ is the Galois closure of a cyclic polynomial $x^4-a$.
\end{itemize}
One easily sees that in the latter case 
any prime $p\equiv 3\mod 4$
will induce a Galois element of order $2$.
So this case does only produce 4-divisible traces at the primes in question
while the first case is ruled out by the assumptions.
This concludes the proof.
\end{proof}

The proposition immediately implies Lemma \ref{lem:4}.
This little extra information will prove crucial for deriving the missing traces:

\subsection{Lemma}
\label{lem:tr}

{\sl
The characteristic polynomials of the Galois representations 
factor as follows at $p=43, 71$:
\begin{eqnarray*}
(V_{15}\otimes W_{15,\ell}(1))
\oplus 
(V_{30}\otimes W_{30,\ell}(1)): &  & \chi(T)= (T^2-d_pT+p)^{27} (T^2+d_pT+p)^{18} \\
V_{45}\otimes W_{45,\ell}(1): &  &  \chi(T)=(T^2-b_pT+p)^{27} (T^2+b_pT+p)^{18} \\
V_{45}'\otimes W_{45,\ell}'(1): &  & \chi(T)= (T^2-c_pT+p)^{27} (T^2+c_pT+p)^{18} 
\end{eqnarray*}
}

\begin{proof}
We return to the concept of Galois traces.
By \ref{ss:rev} the Galois traces of $id, \imath_1, \imath_2, \imath_3$ combine to a $4\times 4$ matrix
\[
M = 
\begin{pmatrix}
u_1 & 9-u_1 & 9 & 9\\
u_2 & 9-u_2 & 9 & 9\\
u_3 & 1-u_3 & 1 & -3\\
u_4 & -1-u_4 & 5 & 1
\end{pmatrix}
\]
Here all $u_i$ range though odd numbers in $\{-15,\hdots,15\}$.
Recall that $M$ can be used to compute a (not necessarily unique) solution to the trace problem 
from point counts over $\FF_p$ and Lefschetz' fixed point formula by \eqref{eq:trg}.
One easily checks that the vector $v_0=(d_p,d_p,b_p,c_p)$ always gives such a solution at $p=43, 71$.
It remains to prove its uniqueness.

The determinant of $M$ is $216(u_1-u_2)$, so uniqueness follows at once if $u_1\neq u_2$.
In the sequel we thus assume $u_1=u_2$, so that $M$ has rank $3$.
The sign argument from \ref{ss:sign} shows that $|u_1|<9$.
This narrows down the possible matrices $M$ enough to compute each one's kernel
which can be used to translate the base solution $v_0$.
Recall that our solution is guaranteed to be integral.
In fact, Lemma \ref{lem:4} leads us to consider a minimal vector $v_\mn$ with entries in $4\ZZ$
spanning the kernel of $M$.
Then we compare its entries against the Weil bounds for the traces of a 2-dimensional Galois representation which (after the Tate twist indicated in the lemma)
give for the translation vector at $p=71$ an absolute upper bound 
\[
|v_\mn|_{\text{max}}\leq 32 =\lfloor 4 \sqrt{71} \rfloor.
\]
This global estimate rules out all but 9 possible translation vectors $v_\mn$.
For each remaining vector, one then checks at $p=43,71$ 
that both
$v_0+v_\mn$ and $v_0-v_\mn$ would have an entry exceeding the Weil bound $2\sqrt p$ for the traces.
Thus $v_0$ is the unique (integral 4-divisible) solution respecting the Weil bounds
regardless of the choice of $M$.
The claimed traces follow.
\end{proof}

\subsection{Conclusion of the proof of Theorem \ref{thmX}}

We have verified that the Galois representations from Theorem \ref{thmX}
have even traces, and in fact the same  traces at any prime from the non-cubic set \ref{ss:Gal}
(as sketched in \ref{ss:aut} for $p\equiv 1\mod 4$, but with the appropriate quadratic twists for $p\equiv 3\mod 4$).
Hence the Galois representations are isomorphic (being semi-simple).
\qed

\subsection{Remark}

In particular the above proof shows that 
the 2-dimensional Galois representations 
$W_{15,\ell,\QQ}$ and $W_{30,\ell,\QQ}$ are isomorphic
as conjectured in \cite[7.3]{BvG}.

\section{Proof of Theorem \ref{thmS}}
\label{s:S}

The proof of Theorem \ref{thmS} proceeds along considerably different lines.
Namely we pursue a genuinely geometric approach in order to
work out a correspondence over $\bar\QQ$ 
between a suitable K3 quotient of $S$ on the one hand
and the product of the  stated elliptic curve and one of its conjugates on the other hand.
In fact, we shall use several isogenous K3 surfaces
following Inose's notion \cite{Inose}
since K3 surfaces connected by a dominant rational map over $\QQ$ share the same
Picard number, $\QQ$-Hodge structure and transcendental Galois representation.

\subsection{K3 quotient}

Maschke's octic $S$ comes with a natural action of the dihedral group $D_4$ of order 8, for instance generated by the involutions
\[
(x_2\mapsto -x_2), \;\; (x_3\mapsto -x_3) \;\; \text{ and } \;\; (x_2,x_3)\mapsto (x_3,x_2).
\]
The quotient surface $S_1$ has a simple model in weighted projective space $\PP[1,1,2,4]$ given by the invariants $x_2^2+x_3^2$ and $x_2^2x_3^2$.
This turns out to be K3, but for now we only need the geometric genus $p_g(S_1)=1$
computed from the invariant 2-forms on $S$.
From the analysis of Maschke's octic we know that the transcendental motive of $S_1$ will be either $W_{3,\ell}', W_{7,\chi_0,\ell}$ or $W_{7,\chi_0,\ell}\otimes\chi_{-1}$ (see \cite[5.7]{BvG}).
Counting points over $\FF_{11}$, for instance, we obtain the trace of F$_{11}^*$ modulo $11$ which agrees exactly with the last Galois representation.
Hence 
\begin{eqnarray}
\label{eq:T}
T(S_1)_\ell = W_{7,\chi_0,\ell}\otimes\chi_{-1}.
\end{eqnarray}
Since all steps in the sequel will be performed over $\QQ$,
the same equality holds for all K3 surfaces derived from $S_1$ in the following.
Throughout we will therefore only refer to the transcendental Galois representation $T_\ell$.

\subsection{Elliptic fibration}

In fact, $S_1$ comes with a natural elliptic fibration where the fibre is given as a double covering of $\PP^1$ ramified in 4 points in terms of the above weighted coordinates.
Completing the square we find a model with a section over $\QQ(\sqrt{-3})$.
Its jacobian $S_2$ over $\QQ$ can be given in the standard coordinates by
\[
S_2:
y^2 = 
x (x^2+21 (t^2+1)^2 x+ 144 (t^4+2 t^3+2 t^2-2 t+1) (t^4-2 t^3+2 t^2+2 t+1)).
\]
This is visibly an elliptic K3 surface with 2-torsion section $(0,0)$
and 8 fibres of Kodaira type $I_2$.
In the following we exhibit various isogenous K3 surfaces
in order to reach one fitting into a Shioda-Inose structure with 
a product of elliptic curves (cf.~\cite{Mo}).
In detail we will derive an isogenous elliptic K3 surface with
singular fibres of  Kodaira type $III^*, I_2, I_1$ twice each and a 2-torsion section
such that the formulas from \cite{GS} will be applicable.

\subsection{Strategy}

A word about the strategy to find a suitable chain of isogenies might be in order.
Generally it seems advisable to reduce the complexity of the elliptic fibrations
such that as much of the N\'eron-Severi group as possible is supported on the singular fibers
(i.e.~the Mordell-Weil rank is lowered).
This not only facilitates explicit calculations,
but also is in accordance with the
low MW-ranks of the known elliptic fibrations involved in Shioda-Inose structures.
In particular this eases the analysis of the configuration of $(-2)$-curves on the K3 surface
in order to single out divisors of the Kodaira types in question.
The linear system of such a divisor 
then induces the desired elliptic fibration.
Working out an explicit  basis of such a linear system
(consisting of a constant and the so-called elliptic parameter $u$)
amounts to linear algebra in the function field of the K3 surface,
but can be tricky nonetheless.
This gives yet another reason to consider auxiliary elliptic fibrations (see \ref{ss:aux}).

In the sequel we first lower the MW-rank to one (to be achieved in \ref{ss:3rd})
and then aim for an elliptic fibration of the specified shape which we will reach in \ref{ss:4th}.
Then we use the Shioda-Inose structure to 
determine the shape of the transcendental Galois representation $T_\ell$ and thus
prove Theorem \ref{thmS}.

\subsection{1st isogeny}

The parameter $u=x/(t^4-2 t^3+2 t^2+2 t+1)$ extracts an alternative elliptic fibration on $S_2$
with $I_0^*$ fibres at $u=0,\infty$.
This fibration comes with 6 further singular fibres of type $I_2$;
its jacobian has full 2-torsion in MW 
and an involution induced from the base by $u\mapsto 12^2/u$.

\subsection{2nd isogeny}

The involution on the base combines with the hyperelliptic involution for a symplectic involution on the jacobian of $S_2$
whose quotient K3 surface $S_3$ comes with the elliptic fibration
\[
S_3: \;\; y^2 = x(x^2+6 (4 t+7) (t-2) (t+2) x+9 (2 t+5) (6 t+13) (t-2)^2 (t+2)^2).
\]
This has singular fibres of type $I_0^*$ at $\infty$ and twice $I_2, I_1^*$ each at
the zeroes of the rightmost polynomial
 as well as
a 2-torsion section $(0,0)$.

\subsection{3rd isogeny}
\label{ss:3rd}

The 2-isogeny results in an elliptic K3 surface $S_4$ with twice $I_1, I_2^*$ 
on top of $I_0^*$ (so that $\rho(S_4)\geq 18$ would follow from the Shioda-Tate formula
and $\rho(S_4)=19$ is equivalent to MW rank one):
\[
S_4:\;\; y^2 = x (x^2 -3 (4 t+7) (t-2) (t+2)x+9 (t-2)^3 (t+2)^3).
\]
Recall that eventually we aim at relating $S$ to an elliptic K3 surface with reducible fibres of type
$III^*, I_2$ twice each.
In fact this can easily be achieved on $S_4$ through an auxiliary fibration
of root type $A_1+E_7+D_8$.
However, it turns out that the final fibration, albeit being jacobian,
does not come with a 2-torsion section, so the results from \cite{GS} do not apply to it.
Instead we shall go through a non-jacobian elliptic fibration on $S_4$ with the right fibre configuration.

\subsection{Auxiliary fibrations}
\label{ss:aux}

In order to extract the above singular fibres on $S_4$ (in a suitable way),
we consider an auxiliary fibration with singular fibres $I_2^*$ twice and $I_2$ four times and full 2-torsion.
This can be given by the parameter $u=x/((t-2)(t+2)^2)$ which results in the Weierstrass form
\[
S_4:\;\; y^2 = x(x-3 t (t-12))(x-t^2 (4t-45)).
\]
The fibre configuration might lead the reader to suspect that this fibration occurs 
on a Kummer surface of product type, confer the classification by Oguiso in \cite{Oguiso}.
However, there is one subtlety:
the torsion sections do not meet the $I_2^*$ fibres in symmetric components.
Hence the fibration cannot be a base change of a rational elliptic surface,
and thus does not coincide with a member of the Kummer fibration.

We continue to extract two singular fibres of type $III^*$ and $I_2$ each.
This is achieved by the elliptic parameter $u=x/(t+90)$
which gives a non-jacobian elliptic fibration on $S_4$.

\subsection{4th isogeny}
\label{ss:4th}

The jacobian of this fibration gives another elliptic K3 surface $S_5$.
Its Weierstrass form can be given as
\[
S_5: \;\; 
y^2 = x(x^2-21 t^2x-t^3(t^2-126t+81)).
\]
A quadratic twist over $\QQ(\sqrt{-21})$ takes this to the form of \cite[(7)]{GS}
with $I_2$ fibres at $a,b=-\frac17\pm\frac 4{49} \sqrt 3$.
Through a Shioda-Inose structure as in \cite[\S 4]{GS} this fibration is related to the product of 
elliptic curves over $K=\QQ(\sqrt 3, \sqrt{-5})$ with j-invariants
\[
j(E) = 192632-111328\sqrt 3-145824\sqrt{-5}+84280\sqrt{-15}
\]
and its conjugate fixing $\sqrt 3$.
One directly verifies that these two curves do not have CM,
but they are 6-isogenous over $\bar\QQ$
(by substituting into the corresponding modular polynomial).
This implies that $\rho(S_5)=\rho(S_1)=19$,
and the transcendental lattice has rank 3 as required.
In fact, the specified elliptic curve $E$ is also 2-isogenous
to its conjugate fixing $\sqrt{-5}$,
and thus 3-isogenous to its conjugate fixing $\sqrt{-15}$.
Comparing j-invariants we find that $E$ is the elliptic $\QQ$-curve $E_{-16/5}$ from \cite[p.~312]{Quer}.

\subsection{Galois representations  of the K3 surfaces}

Since $E$ is $\bar\QQ$-isogenous to its conjugates,
the transcendental motive $T_\ell$ of $S_5$ (and thus of $S_1$)
equals ${\Sym}^2 H^1_\et(E,\QQ_\ell)$ over some finite extension of $\QQ$.
It remains to determine the precise shape of $T_\ell$ over $\QQ$.
The same applies to the K3 surface $X=X_{2/3}$ in the family  $\mathcal X_3$ from \cite{GS}:
the family is parametrised by the modular curve $X^*(6)$, and the given special member
was shown to correspond to $E=E_{-16/5}$ by a similar Shioda--Inose structure.
Thus we can right away compare the 
3-dimensional Galois representations $T_\ell$ of $S_1$ and $T(X)_\ell$ of $X$ over $\QQ$.
Since they agree over some extension while being integral over $\QQ$,
the eigenvalues of Frobenius can only differ by signs.
These signs are encoded in quadratic characters which we compute explicitly below
(higher order characters can only occur in the CM case for the fields $\QQ(\sqrt{-1}), \QQ(\sqrt{-3})$).

\subsection{Determinants}
\label{ss:det}

First we compare the determinants of $T_\ell$ and $T(X)_\ell$.
By \cite[5.7]{BvG} combined with \eqref{eq:T},
$T_\ell$ has determinant
\[
\det(T_\ell) = \chi_5 \chi_\ell^3.
\]
For $X$, we use point counting and Lefschetz' fixed point formula. 
Note that 19 out of 20 eigenvalues $\pm p$ are determined by the 
N\'eron-Severi group of $X$ over $\bar\QQ$ via reduction.
For some primes, the Weil bounds then predict the 20th eigenvalue.
In detail we find
$$
\begin{array}{|c||c|c|c|c|c|c|}
\hline
p & 11 & 53 & 107 & 127 & 139 & 179\\
\hline \text{sign} & - & -  & + & + & - & -\\
\hline
\end{array}
$$
The signs at $11, 53, 107,139$ determine the quadratic character unramified outside $\{2,3,5\}$
uniquely as $\chi_{-5}$.
Since the 20th eigenvalue exactly gives the determinant of $T(X)_\ell$,
we deduce 
\[
\det(T(X)_\ell) = \chi_{-5}\chi_\ell^3.
\]

\subsection{Quadratic twist}
\label{ss:quad}

We are now in the position to compare the signs of the remaining 2 eigenvalues of Frobenius.
For this purpose we define a quadratic character on all primes $p\neq 2,3,5$
which are not supersingular for $E$ (so for instance $p\neq 7$):
\[
\phi: \;\;\; p \mapsto \dfrac{\tr F_p(T_\ell) - \chi_5(p)\, p}{\tr F_p(T(X)_\ell) - \chi_{-5}(p)\,p}.
\]
Numerator and denominator (non-zero by assumption) 
encode the traces of the other 2 eigenvalues, 
thus $\phi$ gives their sign difference.
As a quadratic character unramified outside $\{2,3,5\}$,
$\phi$ is determined by its values at, say,  $\{11,13,17,19\}$.
Counting points we find that $\phi=\chi_{-1}$.
Comparison with the determinants yields that all signs of the eigenvalues
of $T_\ell$ and $T(X)_\ell$ differ by this character:
\begin{eqnarray}
\label{eq:T(X)}
T(X)_\ell = T_\ell \otimes \chi_{-1} = W_{7,\chi_0,\ell}.
\end{eqnarray}

\subsection{Modularity of $E$}
\label{ss:1200}

As the elliptic $\QQ$-curve has an isogeny class defined over $\QQ$,
it comes attached with a Galois representation over $\QQ$.
To work this out,
we fix the Weierstrass model
\begin{eqnarray*}
E:\;\;\; y^2 &  = & {x}^{3}-{\frac {648}{25}}\, \left( 3+\sqrt{-15} \right) ^{2} \left( -11+3\,\sqrt{-15}+8\,\sqrt{-5}-{\frac {24}{5}}\,\sqrt{3} \right) x\\
&&
-{\frac {5184}{125}}\, \left( 3+\sqrt{-15} \right) ^{3}\sqrt{3} \left( 1+\frac 45\,\sqrt{-5} \right)  \left( -{\frac {77}{5}}+9\,\sqrt{-15}+{\frac {56}{5}}\,\sqrt{-5}
\mbox{}-{\frac {72}{5}}\,\sqrt{3} \right).
\end{eqnarray*}
This is the quadratic twist of the curve from \cite[p.~312]{Quer} 
evaluated at $a=-16/5$ over $K\left(\sqrt{\frac{3+\sqrt{-15}}2}\right)$
that makes $E$ 2-isogenous over $K$ to its conjugate $E^\sigma$ 
where $\sigma\in$ Gal$(K/\QQ)$ fixes $\sqrt{-5}$.
Attached to $E$ there is a newform $f{1200} \in S_2(\Gamma_0(1200),\chi_3)$
with character $\chi_3$ and coefficient field $F=\QQ(\sqrt{-2}, \sqrt{-3})$.
Here $f{1200}$ is determined uniquely up to inner twists (accounting for the conjugations in $F$)
by the Hecke eigenvalues
{\renewcommand{\arraystretch}{1.2}
$$
\begin{array}{|r||r|r|r|r|r|r|r|r|r|r|r|r|r|}\hline
p&7&11&13&17&19&23&29&31&37&41&43&47&53 \\
\hline
b_p&
0
& 2\sqrt{6} &
2\sqrt{6}
& 2\sqrt{-3}
& -2\sqrt{-3}
& -6
& -2\sqrt{-2}
& -2\sqrt{-3}
& 2\sqrt{6}
& 4\sqrt{-2}
& 6\sqrt{-2}
& 6
& 6\sqrt{-3}\\
\hline
\end{array}
$$
}
Comparing coefficients and traces at the primes from \ref{ss:Gal},
one directly verifies that
\[
H^1(E,\QQ_\ell) \cong \mbox{Res}^{K}_\QQ V_{f{1200},\ell}
\]
as Gal$(\bar K/K)$-representations.
In particular, this also implies that the symmetric squares of these Galois representations are isomorphic over $K$.

%
%
%

\subsection{Proof of Theorem \ref{thmS}}

We have seen that $T_\ell$ and $\Sym^2 V_{f{1200},\ell}$ agree over some extension.
Hence we can proceed as in \ref{ss:det}, \ref{ss:quad} to determine the precise relation.
By construction, we have
\[
\det {\Sym}^2 V_{f{1200},\ell} = \chi_3 \,\chi_\ell^3.
\]
Then one checks that also the irrational eigenvalues of $T_\ell$ and $\Sym^2 V_{f{1200},\ell}$
 are related by the character $\chi_{15}$
(as they should be since we have two absolutely irreducible Galois representations over $\QQ$).
Thus we deduce
\begin{eqnarray}
\label{eq:T-E}
T_\ell \cong ({\Sym}^2 V_{f{1200},\ell}) \otimes \chi_{15}.
\end{eqnarray}
By \cite[5.7]{BvG} and  \eqref{eq:T},
we know that
\[
T_{S,\ell} = 
(V_{f{15},\ell})^5
\oplus
(T_\ell\otimes\chi_{-1})^{18}
 \oplus
T_\ell^{12}.
\]
Together with 
\eqref{eq:T-E} this combines to the statement of Theorem \ref{thmS}. \qed

\subsection{Remark}
\label{ss:rem}

Note that the above argument proves in particular the modularity of the 
K3 surfaces $S_1,\hdots, S_5$ and $X$ over $\QQ$.
For instance we obtain
\[
T(X)_\ell \cong ({\Sym}^2 V_{f{1200},\ell}) \otimes \chi_{-15}.
\]

{\small

\subsection*{Acknowledgements}

Thanks to Gilberto Bini and Bert van Geemen for the detailed discussions
which lead to this paper,
and to the referee for many helpful comments.
Partial support from ERC through StG 279723 (SURFARI)
is gratefully acknowledged.
}

\end{document}